\documentclass[11pt, reqno]{amsart}
\usepackage{indentfirst, amssymb, amsmath, amsthm, mathrsfs, setspace, indentfirst, enumerate,  mathrsfs, amsmath, amsthm}
\usepackage[bookmarksnumbered, colorlinks, plainpages]{hyperref}
\usepackage{mathrsfs}

\newtheorem*{theoA}{Theorem A}

\newtheorem*{quesA}{Question A}

\newtheorem{ques}{Question}[section]

\newtheorem{theo}{Theorem}[section]
\newtheorem{lem}{Lemma}[section]

\newtheorem{defi}{Definition}[section]

\newcommand{\ol}{\overline}

\newcommand{\be}{\begin{equation}}
\newcommand{\ee}{\end{equation}}
\newcommand{\beas}{\begin{eqnarray*}}
\newcommand{\eeas}{\end{eqnarray*}}
\newcommand{\bea}{\begin{eqnarray}}
\newcommand{\eea}{\end{eqnarray}}

\numberwithin{equation}{section}
\begin{document}

\title[M\MakeLowercase{eromorphic functions together with their partial.....}]{\LARGE M\Large\MakeLowercase{eromorphic functions together with their partial derivatives sharing zeros and poles in} $\mathbb{C}^m$}

\date{}
\author[H. Y. X\MakeLowercase{u}, S. M\MakeLowercase{ajumder} \MakeLowercase{and} G. H\MakeLowercase{aldar}]{H\MakeLowercase{ong} Y\MakeLowercase{an} X\MakeLowercase{u}$^1$, S\MakeLowercase{ujoy} M\MakeLowercase{ajumder}$^2$$^*$ \MakeLowercase{and} G\MakeLowercase{outam} H\MakeLowercase{aldar}$^3$}
\address{$^{1}$ School of Mathematics and Physics, Suqian University, Suqian, Jiangsu 223800, China.}
\email{xhyhhh@126.com}
\address{$^2$$^*$ Department of Mathematics, Raiganj University, Raiganj, West Bengal-733134, India.}
\email{sm05math@gmail.com, sjm@raiganjuniversity.ac.in}
\address{$^3$ Department of Mathematics, Ghani Khan Choudhury Institute of Engineering and Technology, Narayanpur, Malda-732141, West Bengal, India.}
\email{goutamiitm@gmail.com, goutamiit1986@gmail.com}

\renewcommand{\thefootnote}{}
\footnote{*\emph{Corresponding Author}: Sujoy Majumder.}

\renewcommand{\thefootnote}{\arabic{footnote}}
\setcounter{footnote}{0}

\begin{abstract} In the paper, using Nevanlinna's value distribution theory of meromorphic functions in $\mathbb{C}^m$, we study the problem: Let $f$ and $g$ be two finite order non-constant meromorphic functions in $\mathbb{C}^m$ such that $f$ and $g$ share $0$ and $\infty$ CM. If $\partial_{z_i}(f)$ and $\partial_{z_i}(g)$ share $0$ CM and if for one $n>1$, $\partial^n_{z_i}(f)$ and $\partial^n_{z_i}(g)$ share $0$ CM, where $i=1,2,\ldots,m$, then $f$ and $g$ have special presentations. Our finding extends the result of K\"{o}hler (Meromorphic functions sharing zeros and poles and also some of their derivatives sharing zeros, Complex Var., 11 (1989), 39-48) into higher dimensions. Moreover, we show by an example that the restriction on the growth of $f$ and $g$ is necessary.
\end{abstract}

\thanks{Typeset by \AmS -\LaTeX}
\maketitle

\section{{\bf Introduction and main results}}
We define $\mathbb{Z}_+=\mathbb{Z}[0,+\infty)=\{n\in \mathbb{Z}: 0\leq n<+\infty\}$ and $\mathbb{Z}^+=\mathbb{Z}(0,+\infty)=\{n\in \mathbb{Z}: 0<n<+\infty\}$.
On $\mathbb{C}^m$, we define
\[\partial_{z_i}=\frac{\partial}{\partial z_i},\partial^2_{z_jz_i}=\frac{\partial^2}{\partial z_jz_i},\ldots, \partial_{z_i}^{l_i}=\frac{\partial^{l_i}}{\partial z_i^{l_i}}\;\;\text{and}\;\;\partial^{I}=\frac{\partial^{|I|}}{\partial z_1^{i_1}\cdots \partial z_m^{i_m}}\]
where $l_i\in \mathbb{Z}^+\;(i=1,2,\ldots,m)$ and $I=(i_1,\ldots,i_m)\in\mathbb{Z}^m_+$ is a multi-index such that $|I|=\sum_{j=1}^m i_j$.

\medskip
Let $G\not=\varnothing$ be an open subset of $\mathbb{C}^m$. Let $f$ be a holomorphic function on $G$. For a point $a\in\mathbb{C}^m$, we write $f(z)=\sum_{i=0}^{\infty}P_i(z-a)$, where the term $P_i(z-a)$ is either identically zero or a homogeneous polynomial of degree $i$. Denote the zero-multiplicity of $f$ at $a$ by $\mu^0_f(a)=\min\{i:P_i(z-a)\not\equiv 0\}$ in the sense of \cite[Definition 2.1]{HF}.

\medskip
Let $f$ be a meromorphic function on $G$. Then there exist holomorphic functions $g$ and $h$ such that $hf=g$ on $G$ and $\dim_z I_f\leq m-2$, where $I_f=\{z\in \mathbb{C}^m: g(z)=h(z)=0\}$. Therefore the $c$-multiplicity of $f$ is just $\mu^c_f=\mu^0_{g-ch}$ if $c\in\mathbb{C}$ and $\mu^c_f=\mu^0_h$ if $c=\infty$. The function $\mu^c_f: \mathbb{C}^m\to \mathbb{Z}$ is nonnegative and is called the $c$-divisor of $f$. If $f\not\equiv 0$ on each component of $G$, then $\nu=\mu_f=\mu^0_f-\mu^{\infty}_f$ is called the divisor of $f$.
For $k\in\mathbb{N}$, define the truncated multiplicity functions on $\mathbb{C}^m$ by $\mu_{f,k}^a(z)=\min\{\mu_f^a(z),k\}$.

\smallskip
Let $f$, $g$ and $a$ be meromorphic functions on $\mathbb{C}^m$. Then one can find three pairs
of entire functions $f_1$ and $f_2$, $g_1$ and $g_2$, and $a_1$ and $a_2$, in which each pair is coprime
at each point in $\mathbb{C}^m$ such that $f = f_2/f_1$, $g=g_2/g_1$ and $a = a_2/a_1$.
We say that $f$ and $g$ share $a$ by counting multiplicities (CM) if
\[\mu_{a_1f_2-a_2f_1}^0=\mu_{a_1g_2-a_2g_1}^0\;(a\not\equiv \infty)\;\;\text{and}\;\;\mu_{f_1}^0=\mu_{g_1}^0\;\;(a=\infty).\]

Again we say that $f$ and $g$ share $a$ by ignoring multiplicities (IM) if
\[\mu_{a_1f_2-a_2f_1,1}^0=\mu_{a_1g_2-a_2g_1,1}^0\;(a\not\equiv \infty)\;\;\text{and}\;\;\mu_{f_1,1}^0=\mu_{g_1,1}^0\;\;(a=\infty).\]

\medskip
Understanding the relationship between meromorphic functions by their shared zeros and poles is a key area of research in complex analysis. This allows for the characterization of functions based on their value distribution, particularly in relation to their derivatives, a field known as Nevanlinna Theory. This helps in problems such as determining if a function is uniquely determined by its zeros and poles.

\medskip
Nevanlinna theory gives powerful tools to compare the growth and value distribution of meromorphic functions. A foundational result is that if two meromorphic functions share five distinct values IM, they must be identical.

Two meromorphic functions $f(z)$ and $g(z)$ in $\mathbb{C}$ sharing the same zeros and poles with their multiplicities implies they differ only by a constant multiplier. More generally, if their zeros and poles are the CM shared, they are closely related through a limited set of transformations, such as being related by a M\"{o}bius transformation or an exponential function. Research in this area uses Nevanlinna theory to understand how sharing zeros and poles with derivatives can determine the relationship between functions. In the paper, we are interested in the theory behind meromorphic functions sharing zeros and poles, and also the scenario where some derivatives of these functions share zeros. This falls under the broader framework of uniqueness theory, particularly using Nevanlinna theory and value distribution.

\smallskip
Research on ``Meromorphic functions sharing zeros and poles and also some of their derivatives sharing zeros" focuses on the uniqueness of meromorphic functions by establishing conditions under which two non-constant meromorphic functions $f(z)$ and $g(z)$ in $\mathbb{C}$, must be related, when they share specific zeros and poles, and when their derivatives share zeros. Researchers study conditions under which $f(z)$ and $g(z)$ are related if their derivatives $f^{(k)}(z)$ and $g^{(k)}(z)$ share certain values CM or IM. By weakening the sharing conditions for the functions and their derivatives, mathematicians can derive more precise and comprehensive conclusions about the relationship between these functions. This area has applications in complex dynamics:\par
$(a)$ Studying how iterations of meromorphic functions behave,\par
$(b)$ Functions sharing poles/zeros with their iterates can exhibit special types of Julia sets or Fatou sets (stable/chaotic regions)\\
and in differential equations:\par
$(a)$ If solutions to a differential equation share poles/zeros or those of their derivatives, that can help classify all solutions or establish uniqueness,\par
$(b)$ Used in Painlev$\acute{e}$ equations, where meromorphic solutions have special value-sharing properties.

\medskip
For non-constant entire functions $f(z)$ and $g(z)$ in $\mathbb{C}$ with finite order, C. C. Yang \cite{CCY} and G. G. Gundersen \cite{GGG} studied the case where $f^{(j)}(z)$ and $g^{(j)}(z)$ share 0 CM for $j=0,1$ and examined what implications this strong form of sharing (CM sharing of $0$ for both the functions and their first derivatives) has on the relationship between $f(z)$ and $g(z)$.

\medskip
For non-constant meromorphic functions $f(z)$ and $g(z)$ in $\mathbb{C}$, we know that $f^{(j)}(z)$ and $g^{(j)}(z)$ share $0$ and $\infty$ CM for each non-negative integer $j$ whenever $f(z)$ and $g(z)$ satisfy one of the following four conditions:
\begin{enumerate}
\item[(i)] $f=cg$, where $c\in\mathbb{C}\backslash \{0\}$;

\smallskip
\item[(ii)] $f(z)=e^{az+b}$ and $g(z)=e^{cz+d}$, where $a,c\in\mathbb{C}\backslash \{0\}$ and $b,d\in\mathbb{C}$;

\smallskip
\item[(iii)] $f(z)=a\left(1-be^{cz}\right)$ and $g(z)=d\left(e^{-cz}-b\right)$, where $a,b,c,d\in\mathbb{C}\backslash \{0\}$;

\smallskip
\item[(iv)] $f(z)=a/\left(1-be^{\beta(z)}\right)$ and $g(z)=a/\left(e^{-\beta(z)}-b\right)$, where $a,b\in\mathbb{C}\backslash \{0\}$ and $\beta$ is a non-constant entire function.
\end{enumerate}

\smallskip
In 1984, A. Hinkkanen \cite[p. 492]{BBH} proposed the following problem:
\begin{quesA} {\bf (Hinkkanen's Problem)} Does there exist an integer $n\geq 2$ such that $f(z)$ and $g(z)$ satisfy one of the conditions $(i)$-$(iv)$ when $f^{(j)}(z)$ and $g^{(j)}(z)$ share the values $0$ and $\infty$ CM for $j=0,1,\ldots,n$?
\end{quesA}

Hinkkanen's Problem refers to a notable result in the uniqueness theory of meromorphic functions. Hinkkanen originally raised this question in the context of how many derivative zero sets are needed to uniquely determine a meromorphic function.

\smallskip
In 1989, K\"{o}hler \cite{LK} proved that $n=6$ solves the problem. In 1990, Tohge \cite{KT} considered the case $n=2,3$ under restrictions on the growth of $f(z)$ and $g(z)$. According to the answer to Hinkkanen's Problem, Frank et al. \cite{FHV} demonstrated that a meromorphic function may be uniquely identified if we know the positions of its zeros, poles, and the zeros of its first four derivatives. In fact, Frank et al. \cite{FHV} proved that: The sharp answer to Hinkkanen's problem is $n=4$.
In Example 1 \cite{FHV}, Frank et al. \cite{FHV} demonstrated the sharpness of Hinkkanen's Problem for $n=4$.

\smallskip
In 1989, K\"{o}hler \cite{LK} examined Hinkkanen's Problem with the additional condition that the order of $f(z)$ and $g(z)$ are finite. He proved the following result.

\begin{theoA}\cite[Theorem 2]{LK} Let $f(z)$ and $g(z)$ be two non-constant meromorphic functions in $\mathbb{C}$ of finite order such that $f(z)$ and $g(z)$ share $0$ and $\infty$ CM. Suppose $f^{(1)}(z)$ and $g^{(1)}(z)$ share $0$ CM and for $k>1$, $f^{(k)}(z)$ and $g^{(k)}(z)$ share $0$ CM. Then $f(z)$ and $g(z)$ satisfy one of the above case $(i)$-$(iv)$ and $\beta(z)$ in $(iv)$ is a non-constant polynomial.
\end{theoA}

By the example $f(z)=\exp(\exp (z))$ and $g(z)=\exp(\exp(-z))$, it is clear that Theorem A does not hold for functions of infinite order.

\medskip
Theorem A leads to the following natural question:
\begin{ques} Does Theorem A work in several complex variables?
\end{ques}

We respond to Question 1.1 in the paper with the following outcome.

\begin{theo}\label{t2.1} Let $f$ and $g$ be two non-constant finite order meromorphic function in $\mathbb{C}^m$ such that $f$ and $g$ share $0$ and $\infty$ CM. Suppose $\partial_{z_i}(f)$ and $\partial_{z_i}(g)$ share $0$ CM and for one $n>1$, $\partial_{z_i}^n(f)$ and $\partial_{z_i}^n(g)$ share $0$ CM, where $i=1,2,\ldots,m$. Then $f(z)$ and $g(z)$ satisfy one of the following four conditions:
\begin{enumerate}
\item[(1)] $f=cg$, where $c\in\mathbb{C}\backslash \{0\}$;

\smallskip
\item[(ii)] $f(z)=e^{a_1z_1+a_2z_2+\ldots+a_m z_m+a}$ and $g(z)=e^{b_1z_1+b_2z_2+\ldots+b_mz_m+b}$, where $a,a_1,\ldots,a_m,b_1,\ldots,\\b_m,b\in\mathbb{C}$ such that $(a_1,\ldots,a_m)\neq (0,\ldots,0)$ and $(b_1,\ldots,b_m)\neq (0,\ldots,0)$.

\smallskip
\item[(iii)]$f(z)=a\left(1-be^{(a_1z_1+a_2z_2+\ldots+a_mz_m)}\right)$ and $g(z)=c\left(e^{-(a_1z_1+a_2z_2+\ldots+a_mz_m)}-b\right)$,
where $a,b,c\in\mathbb{C}\backslash \{0\}$ and $a_1,a_2,\ldots,a_m\in\mathbb{C}$ such that $(a_1,a_2,\ldots,a_m)\neq (0,0,\ldots,0)$;

\smallskip
\item[(iv)] $f(z)=\frac{d}{1-ce^{-P(z)}}$ and $g(z)=\frac{d}{e^{P(z)}-c}$, where $c,d\in\mathbb{C}\backslash \{0\}$ and $P$ is a non-constant polynomial in $\mathbb{C}^m$.
\end{enumerate}
\end{theo}

Let
\[f(z)=\left(e^{2(z_1+\ldots+z_m)}-1\right)\exp\left(-ie^{z_1+\ldots+z_m}\right)\]
and
\[g(z)=\left(1-e^{-2(z_1+\ldots+z_m)}\right)\exp\left(ie^{-(z_1+\ldots+z_m)}\right).\]

Clearly $f$ and $g$ share $0$ and $\infty$ CM, $\partial_{z_i}(f)$ and $\partial_{z_i}(g)$ share $0$ CM and $\partial^2_{z_i}(f)$ and $\partial^2_{z_i}(g)$ share $0$, where $i=1,2,\ldots,m$, but $f(z)$ and $g(z)$ do not satisfy any case of $(i)$-$(iv)$ of Theorem \ref{t2.1}. Consequently for meromorphic functions of infinite order Theorem \ref{t2.1} is invalid.

\medskip
Before going to the proof of our main result, we will briefly recall some of the standard notation of Nevanlinna theory in $\mathbb{C}^m$ (see \cite{HLY1,ST1}).
Nevanlinna theory, originally formulated for meromorphic functions of a single complex variable, has been extended to the setting of several complex variables to study the value distribution of meromorphic maps from $\mathbb{C}^m$ into complex projective spaces or complex manifolds. It provides profound insights into: $(a)$ The growth and value distribution of meromorphic mappings; $(b)$ Uniqueness and rigidity phenomena for solutions of linear and nonlinear PDEs and $(c)$ Criteria for compactness and convergence in families of meromorphic maps. Its applications span complex analysis, partial differential equations, complex geometry, and mathematical physics, offering both powerful techniques and elegant results. These references \cite{Banerjee-Majumder-2026}, \cite{CK1}, \cite{PVD1}-\cite{PVD2}, \cite{hy}, \cite{BQL1}-\cite{FL1}, \cite{Majumder-2026}, \cite{Majumder-Das}, \cite{Majumder-Das-Pramanik-2025}, \cite{Majumder-Sarkar}, \cite{Majumder-Sarkar-2026}, \cite{Majumder-Sarkar-2027}, \cite{MSP}, \cite{GS1}, \cite{GS2}, \cite{XLL1}, \cite{XW1} and \cite{XH}
provide a foundation for understanding the current state of research in Nevanlinna value distribution theory in several complex variables.

\smallskip
On $\mathbb{C}^m$, the exterior derivative $d$ splits $d= \partial+ \bar{\partial}$ and twists to $d^c= \frac{\iota}{4\pi}\left(\bar{\partial}- \partial\right)$. Clearly $dd^{c}= \frac{\iota}{2\pi}\partial\bar{\partial}$. A non-negative function $\tau: \mathbb{C}^m\to \mathbb{R}[0,b)\;(0<b\leq \infty)$ of class $\mathbb{C}^{\infty}$ is said to be an exhaustion of $\mathbb{C}^m$ if $\tau^{-1}(K)$ is compact whenever $K$ is.
An exhaustion $\tau_m$ of $\mathbb{C}^m$ is defined by $\tau_m(z)=||z||^2$. The standard Kaehler metric on $\mathbb{C}^m$ is given by $\upsilon_m=dd^c\tau_m>0$. On $\mathbb{C}^m\backslash \{0\}$, we define
\[\omega_m=dd^c\log \tau_m\geq 0\;\;\text{and}\;\;\sigma_m=d^c\log \tau_m \wedge \omega_m^{m-1}.\]

For any $S\subseteq \mathbb{C}^m$, let $S[r]$ and $S\langle r\rangle$ be the intersection of $S$ with respectively the closed ball and the sphere of radius $r>0$ centered at $0\in \mathbb{C}^m$. Also for $t>0$, the counting function $n_{\nu}$ is defined by
\beas n_{\nu}(t)=t^{-2(m-1)}\int_{A[t]}\nu \upsilon_m^{m-1},\eeas
where $A=\text{supp}\;\nu=\text{supp}\;\mu_f=\ol{\{z\in G: \nu(z)\neq 0\}}$.
The valence function of $\nu$ is defined by
\[N_{\nu}(r)=N_{\nu}(r,r_0)=\int_{r_0}^r n_{\nu}(t)\frac{dt}{t}\;\;(r\geq r_0).\]

\medskip
For $a\in\mathbb{P}^1=\mathbb{C}\cup\{\infty\}$, we write $N_{\mu_f^a}(r)=N(r,a;f)$ if $a\in\mathbb{C}$ and $N_{\mu_f^a}(r)=N(r,f)$ if $a=\infty$.
For $k\in\mathbb{N}$, define the truncated multiplicity functions on $\mathbb{C}^m$ by
\beas \bar{\mu}_{f=k}^a(z)=\begin{cases}
1, &\text{if $\mu_f^a(z)=k$}\\
0, &\text{if $\mu_f^a(z)\neq k$}
\end{cases}\;\;\text{and}\;\;
\bar{\mu}_{f(k}^a(z)=\begin{cases}
1, &\text{if $\mu_f^a(z)\geq  k$}\\
0, &\text{if $\mu_f^a(z)<k$}
\end{cases}\eeas
and the truncated valence functions
\beas N_{\nu}(t)=\begin{cases}
\ol{N}(t,a;f), &\text{if $\nu=\mu_{f,1}^a$}\\
\ol{N}_{(k}(t,a;f), &\text{if $\nu=\bar{\mu}_{f(k}^a$}.
\end{cases}\eeas

\begin{defi} An algebraic subset $X$ of $\mathbb{C}^m$ is defined as a subset
\[X=\left\lbrace z\in\mathbb{C}^m: P_j(z)=0,\;1\leq j\leq l\right\rbrace\]
with finitely many polynomials $P_1(z),\ldots, P_l(z)$.
A divisor $\nu$ on $\mathbb{C}^m$ is said to be algebraic if $\nu$ is the zero divisor of a polynomial. In this case the counting
function $n_{\nu}$ is bounded (see \cite{GK1,ST1}).
\end{defi}

\smallskip
With the help of the positive logarithm function, we define the proximity function of $f$ by
\[m(r, f)=\mathbb{C}^m\langle r; \log^+ |f| \rangle=\int_{\mathbb{C}^m\langle r\rangle} \log^+ |f|\;\sigma_m.\]

The characteristic function of $f$ is defined by $T(r, f)=m(r,f)+N(r,f)$.
We define $m(r,a;f)=m(r,f)$ if $a=\infty$ and $m(r,a;f)=m(r,1/(f-a))$ if $a$ is finite complex number. If $a\in\mathbb{C}$, then the first main theorem states that $m(r,a;f)+N(r,a;f)=T(r,f)+O(1)$, where $O(1)$ denotes a bounded function when $r$ is sufficiently large. We define the order of $f$ by
\[\rho(f):=\limsup _{r \rightarrow \infty} \frac{\log T(r, f)}{\log r}.\]

\begin{defi} Let $U\subset \mathbb{C}^m$ be an open subset. A closed subset $A\subset U$ is said to be analytic,
if for each $a\in A$, there are a finite number of holomorphic functions $f_1(z),f_2(z),\ldots,f_l(z)$ defined in a neighbourhood $N(a)$ of $a$ such that (see \cite[pp. 42]{NW})
\[A\cap N(a)=\{z\in N(a):f_1(z)=f_2(z)=\cdots=f_l(z)=0\}.\]
If $f_j(z),1\leq j\leq l$ can be taken so that their differentials at $a\in A$, $df_1(a), df_2(a), \ldots,df_l(a)$ are linearly independent, then $a$ is called a regular point of $A$. A point of $A$ which is not regular called a singular point.
\end{defi}

\section {{\bf Auxiliary lemmas}}

In the proof of Theorem \ref{t2.1}, we use the following key lemmas.

\begin{lem}\label{l2.1}\cite[Lemma 1.37]{HLY1} Let $f:\mathbb{C}^m\to\mathbb{P}^1$ be a non-constant meromorphic function and let $I=(i_1,\ldots,i_m)\in \mathbb{Z}^m_+$ be a multi-index. Then for any $\varepsilon>0$, we have
\[m\left(r,\frac{\partial^I(f)}{f}\right)\leq |I|\log^+T(r,f)+|I|(1+\varepsilon)\log^+\log T(r,f)+O(1),\]
for all large $r$ outside a set $E$ with $\text{log mes}\;E=\int_E d\log r<\infty$.
\end{lem}

\begin{lem}\label{L.0}\cite[Proposition 3.2]{HY1} Let $P$ be a non-constant entire function in $\mathbb{C}^m$. Then
\[\rho(e^P)=
\begin{cases}
\deg(P) & \text{if $P$ is a polynomial,}\\
+\infty & \text{otherwise}.
\end{cases}\]
\end{lem}
\begin{lem}\label{L.1} Let $F$ be a non-constant meromorphic function in $\mathbb{C}^m$ and $f_k=\frac{\partial_{z_k}(F)}{F}$, where $k=1,2,\ldots,m$. Then
\beas \frac{\partial_{z_k}^l(F)}{F}=f_k^{l}+\frac{l(l-1)}{2}f_k^{l-2}\partial_{z_k}(f_k)+\tilde A_lf_k^{l-3}\partial^2_{z_k}(f_k)+\tilde B_lf_k^{l-4}(\partial_{z_k}(f_k))^2+P_{l-3}(f_k),\eeas
where
\[\tilde A_l=\frac{l(l-1)(l-2)}{6},\;\;\tilde B_l=\frac{l(l-1)(l-2)(l-3)}{8},\]
and $P_{l-3}(f_k)$ is a differential polynomial with constant coefficients, which vanishes identically for $l\leq 3$ and has degree $l-3$ when $l>3$.
\end{lem}
\begin{proof} We have
\bea\label{eL.1}f_k=\frac{\partial_{z_k}(F)}{F},\eea
where $k=1,2,\ldots,m$. Clearly from (\ref{eL.1}), we have
\beas \partial_{z_k}^2(F)=f_k\partial_{z_k}(F)+F\partial_{z_k}(f_k)&=&f_k^2F+F\partial_{z_k}(f_k)
=(f_k^2+\partial_{z_k}(f_k))F,\eeas
i.e,
\[\frac{\partial_{z_k}^2(F)}{F}=f_k^2+\partial_{z_k}(f_k).\]

We inductively deduce (see \cite[Lemma 3.5]{wkh2}) that
\beas \frac{\partial_{z_k}^l(F)}{F}=f_k^{l}+\frac{l(l-1)}{2}f_k^{l-2}\partial_{z_k}(f_k)+\tilde A_lf_k^{l-3}\partial^2_{z_k}(f_k)+\tilde B_lf_k^{l-4}(\partial_{z_k}(f_k))^2+P_{l-3}(f_k),\eeas
where
\[\tilde A_l=\frac{l(l-1)(l-2)}{6},\;\;\tilde B_l=\frac{l(l-1)(l-2)(l-3)}{8},\]
and $P_{l-3}(f_k)$ is a differential polynomial with constant coefficients, which vanishes identically for $l\leq 3$ and has degree $l-3$ when $l>3$.
\end{proof}

\begin{lem}\label{L.2} Let $f$ and $g$ be two non-constant finite order meromorphic function in $\mathbb{C}^m$ such that $f$ and $g$ share $0$ and $\infty$ CM. Suppose $\partial_{z_i}(f)$ and $\partial_{z_i}(g)$ share $0$ CM and for one $n>1$, $\partial_{z_i}^n(f)$ and $\partial_{z_i}^n(g)$ share $0$ CM, where $i=1,2,\ldots,m$. If either $\partial^n_{z_i}(f)\equiv 0$ or $\partial_{z_i}(f)/f$ is a rational function for all $i=1,2,\ldots,m$, then we have either $f\equiv cg$, where $c\in\mathbb{C}\backslash \{0\}$ or
\[f(z)=e^{a_1z_1+a_2z_2+\ldots+a_m z_m+a}\;\;\text{and}\;\;g(z)=e^{b_1z_1+b_2z_2+\ldots+b_mz_m+b},\]
where $a,a_1,\ldots,a_m,b_1,\ldots,b_m,b\in\mathbb{C}$ such that $(a_1,\ldots,a_m)\neq (0,\ldots,0)$ and $(b_1,\ldots,b_m)\neq (0,\ldots,0)$.
\end{lem}
\begin{proof} By the given conditions, we have $f$ and $g$ share $0$ and $\infty$ CM, $\partial_{z_i}(f)$ and $\partial_{z_i}(g)$ share $0$ CM and $\partial_{z_i}^n(f)$ and $\partial_{z_i}^n(g)$ share $0$ CM, where $n>1$ and $i=1,2,\ldots,m$.

Now we consider following two cases.\par

\medskip
{\bf Case 1.} Suppose $\partial^n_{z_i}(f(z))\equiv 0$ for all $i=1,2,\ldots,m$. Since $\partial_{z_i}^n(f)$ and $\partial_{z_i}^n(g)$ share $0$ CM, it follows that $\partial^n_{z_i}(g(z))\equiv 0$ for all $i=1,2,\ldots,m$. Therefore we have
\bea\label{n.1} \frac{\partial^nf(z)}{\partial z_i^n}\equiv 0\;\;\text{and}\;\;\frac{\partial^ng(z)}{\partial z_i^n}\equiv 0,\eea
for $i=1,2,\ldots,m$.

We now want to prove that $f(z)$ is a polynomial in $\mathbb{C}^m$. Here we prove that $f(z)$ is a polynomial in $\mathbb{C}^m$ by induction on the number $m$ of variables. If $m=1$, from (\ref{n.1}), we get that $f(z)$ is a polynomial in $\mathbb{C}$. Let us suppose $m=2$. Since $\frac{\partial^nf(z_1,z_2)}{\partial z_2^n}\equiv 0$, on integration, we have
\bea\label {n.2} f(z_1,z_2)=\phi_1(z_1)z_2^{n-1}+\phi_2(z_1)z_2^{n-2}+\cdots+\phi_n(z_1),\eea
where $\phi_i(z_1)$'s are entire functions in $\mathbb{C}$ in the variable $z_1$. Note that $\frac{\partial^n f(z_1,z_2)}{\partial z_1^n}\equiv 0$ and so (\ref{n.2}) gives
\[\phi_1^{(n)}(z_1)z_2^{n-1}+\phi_2^{(n)}(z_1)z_2^{n-2}+\cdots+\phi_n^{(n)}(z_1)\equiv 0,\]
which shows that $\phi_i^{(n)}(z_1)\equiv 0$ for $i=1,2,\dots,n$. Therefore on integration, we get
\[\phi_i(z_1)=c_{n-1,i}z_1^{n-1}+c_{n-2,i}z_1^{n-2}+c_{0,i},\]
where $c_{i,j}\;(0\leq i\leq n-1, j=1,2,\ldots,n)$ are constants in $\mathbb{C}$. Then from (\ref{n.2}), we get
\[f(z_1,z_2)=(c_{n-1,1}z_1^{n-1}+c_{n-2,1}z_1^{n-2}+c_{0,1})z_2^{n-1}+\cdots+(c_{n-1,n}z_1^{n-1}+c_{n-2,n}z_1^{n-2}+c_{0,n}),\]
 which shows that $f(z_1,z_2)$ is a polynomial in $\mathbb{C}^2$. Now we fix $m\geq 2$ and assume that $f(z)$ is a polynomial for variables of number at most $m-1$. Since $\frac{\partial^nf(z)}{\partial z_m^n}\equiv 0$, we have
\bea\label{n.3} f(z_1,z_2,\ldots,z_m)&=& A_1(z_1,z_2,\ldots,z_{m-1})z_m^{n-1}+A_2(z_1,z_2,\ldots,z_{m-1})z_m^{n-2}+\cdots\\&&+A_n(z_1,z_2,\ldots,z_{m-1}).\nonumber\eea

Now using (\ref{n.1}) to (\ref{n.3}), we get
\bea\label{n.4} \frac{\partial^n A_1(z_1,\ldots,z_{m-1})}{\partial z_i^n}z_m^{n-1}+\frac{\partial^2 A_2(z_1,\ldots,z_{m-1})}{\partial z_i^2}z_m^{n-2}+\cdots+\frac{\partial^n A_n(z_1,\ldots,z_{m-1})}{\partial z_i^n}\equiv 0\eea
for $i=1,2,\ldots,m-1$. Therefore (\ref{n.4}) yields
\bea\label{n.5} \frac{\partial^n A_j(z_1,z_2,\ldots,z_{m-1})}{\partial z_i^n}\equiv 0\eea
for all $j=1,2,\ldots,n$, where $i=1,2,\ldots,m-1$. Then by the induction assumptions, we find from (\ref{n.5}) that  $A_i(z_1,z_2,\ldots,z_{m-1})$ $(\forall\; i=1,2,\ldots,n)$ are polynomials in the variables $z_1,z_2,\ldots,z_m$. Therefore from (\ref{n.3}), we see that $f(z_1,z_2,\ldots,z_m)$ is a polynomial in $z_1,z_2,\ldots,z_m$. Similarly $g(z_1,z_2,\ldots,z_m)$ is a polynomial in $z_1,z_2,\ldots,z_m$. Since $f$ and $g$ share $0$ CM, it follows that $f=cg$, where $c\in\mathbb{C}\backslash \{0\}$.\par

\medskip
{\bf Case 2.} Suppose $\partial_{z_i}(f)/f$ is a rational function for all $i=1,2,\ldots,m$. Since $f$ and $g$ share $0$ and $\infty$ CM, there exists a polynomial $P$ in $\mathbb{C}^m$ such that $f=ge^P$. Clearly, we have
\[\frac{\partial_{z_i}(f)}{f}=\frac{\partial_{z_i}(g)}{g}+\partial_{z_i}(P),\]
which shows that $\partial_{z_i}(g)/g$ is also a rational function for all $i=1,2,\ldots,m$. By the given conditions, we have $f$ and $g$ share $0$ and $\infty$ CM and $\partial_{z_i}(f)$ and $\partial_{z_i}(g)$ share $0$ CM for all $i=1,2,\ldots,m$. Consequently $\partial_{z_i}(f)/f$ and $\partial_{z_i}(g)/g$ share $0$ and $\infty$ CM. As a result, we have
\bea\label{n.5a}\frac{\partial_{z_i}(f)}{f}=A_i\frac{\partial_{z_i}(g)}{g},\eea
where $A_i\in\mathbb{C}\backslash \{0\}$ and $i=1,2,\ldots,m$. On integration, we get
\bea\label{n.6} f(z)=g^{A_i}(z)e^{g_i(z)},\eea
where $g_i(z)=g(z_1,z_2,\ldots,z_{i-1},z_{i+1},\ldots,z_m)$ is a polynomial for $i=1,2,\ldots,m$.

Now we consider following two sub-cases.\par

\medskip
{\bf Sub-case 2.1.} Suppose $0$ and $\infty$ are not both Picard exceptional values of $f$. Since $f$ and $g$ share $0$ and $\infty$ CM, it follows that $0$ and $\infty$ are not both Picard exceptional values of $g$. Now from (\ref{n.6}), we deduce that $A_i=1$, where $i=1,2,\ldots,m$ and so from (\ref{n.5}) and (\ref{n.6}), we get respectively
\bea\label{n.7aa}\frac{\partial_{z_j}(f)}{f}=\frac{\partial_{z_j}(g)}{g},\eea
where $j=1,2,\ldots,m$ and
\bea\label{n.7} f(z)=g(z)e^{g_i(z)},\eea
where $i=1,2,\ldots,m$. Taking partial derivative with respect to $z_j$ on the both sides of (\ref{n.7}) and using (\ref{n.7aa}), we get
$\partial_{z_j}(g_i(z))\equiv 0$, where $i,j=1,2,\ldots,m$ such that $i\neq j$. Since $g_i(z)$ is independent of $z_i$, we conclude that $g_i(z)$ is a constant. Now in view of (\ref{n.7}), we can assume that $e^{g_i(z)}=c$ for all $i=1,2,\ldots,m$. Therefore $f=cg$, where $c\in\mathbb{C}\backslash \{0\}$.\par

\medskip
{\bf Case 2.2.} Suppose $0$ and $\infty$ are both Picard exceptional values of $f$. Since $f$ and $g$ share $0$ and $\infty$ CM, it follows that $0$ and $\infty$ are both Picard exceptional values of $g$. Consequently $\partial_{z_i}(f)/f$ and $\partial_{z_i}(g)/g$ are polynomials in $\mathbb{C}^m$ for all $i=1,2,\ldots,m$.
Let $g=\exp(Q)$, where $Q$ is a polynomial in $\mathbb{C}^m$. Clearly
\bea\label{n.7a} \frac{\partial_{z_i}(g)}{g}=\partial_{z_i}(Q).\eea

Now from (\ref{n.5a}), we get
\bea\label{n.7b}\frac{\partial_{z_i}(f)}{f}=A_i\partial_{z_i}(Q)\eea
and so on integration, we have
\bea\label{n.8} f(z)=e^{A_iQ(z)+H_i(z)},\eea
where $H_i(z)=H(z_1,z_2,\ldots,z_{i-1},z_{i+1},\ldots,z_m)$ is a polynomial, $i=1,2,\ldots,m$. Note that (\ref{n.8}) gives
\bea\label{n.8a}\frac{\partial_{z_j}(f)}{f}=A_i\partial_{z_j}(Q)+\partial_{z_j}(H_i),\eea
where $i,j\in\{1,2,\ldots,m\}$ such that $i\neq j$.

\smallskip
By Lemma \ref{L.1}, we can say that $\partial_{z_i}^n(f)/f$ and $\partial_{z_i}^n(g)/g$ are rational functions for all $i=1,2,\ldots,m$. Also by the given conditions, we have $f$ and $g$ share $0$ and $\infty$ CM and $\partial_{z_i}^n(f)$ and $\partial_{z_i}^n(g)$ share $0$ CM for all $i=1,2,\ldots,m$. Consequently $\partial_{z_i}^n(f)/f$ and $\partial_{z_i}^n(g)/g$ share $0$ and $\infty$ CM. As a result, we have
\bea\label{n.9}\frac{\partial_{z_i}^n(f)}{f}=B_i\frac{\partial_{z_i}^n(g)}{g},\eea
where $B_i\in\mathbb{C}\backslash \{0\}$, where $i=1,2,\ldots,m$. Now in view of (\ref{n.7a}) and (\ref{n.7b}) and using Lemma \ref{L.1} to (\ref{n.9}), we have
\bea\label{n.10} &&(A_i\partial_{z_i}(Q))^n+\frac{n(n-1)}{2}(A_i\partial_{z_i}(Q))^{n-2}\partial_{z_i}(A_i\partial_{z_i}(Q))+\tilde A_n(A_i\partial_{z_i}(Q))^{n-3}\partial^2_{z_i}(A_i\partial_{z_i}(Q))\nonumber\\&&+\tilde B_n(A_i\partial_{z_i}(Q))^{n-4}(\partial_{z_i}(A_i\partial_{z_i}(Q)))^2+P_{n-3}(A_i\partial_{z_i}(Q))\\&&=
B_i\Big(\partial_{z_i}(Q))^n+\frac{n(n-1)}{2}(\partial_{z_i}(Q))^{n-2}\partial_{z_i}(\partial_{z_i}(Q))+\tilde A_n(\partial_{z_i}(Q))^{n-3}\partial^2_{z_i}(\partial_{z_i}(Q))\nonumber\\ &&+\tilde B_n(\partial_{z_i}(Q))^{n-4}(\partial_{z_i}(\partial_{z_i}(Q)))^2+P_{n-3}(\partial_{z_i}(Q))\Big)\nonumber
\eea
where $i=1,2,\ldots,m$,
\[\tilde A_n=\frac{n(n-1)(n-2)}{6},\;\;\tilde B_n=\frac{n(n-1)(n-2)(n-3)}{8},\]
and $P_{n-3}(\partial_{z_i}(P))$ is a differential polynomial with constant coefficients, which vanishes identically for $n\leq 3$ and has degree $n-3$ when $n>3$.

We consider following two sub-cases.\par

\medskip
{\bf Sub-case 2.2.1.} Suppose $\deg(\partial_{z_i}(Q))\geq 1$ for all $i\in\{1,2,\ldots,m\}$. Now from (\ref{n.10}), we deduce that $A_i^n=B_i$ and $A_i^{n-1}=B_i$ for all $i\in\{1,2,\ldots,m\}$. Clearly $A_i=1$ for all $i\in\{1,2,\ldots,m\}$.
Now proceeding in the same way as done in the proof of Sub-case 2.1, we easily get $f=cg$, where $c\in\mathbb{C}\backslash \{0\}$.\par

\medskip
{\bf Sub-case 2.2.2.} Suppose $\deg(\partial_{z_i}(Q))<1$ for at least one $i\in\{1,2,\ldots,m\}$.
We know that $\partial_{z_j}(f)/f$ and $\partial_{j}(g)/g$ share $0$ CM, where $i,j\in\{1,2,\ldots,m\}$ such that $i\neq j$.
Therefore from (\ref{n.7a}) and (\ref{n.8a}), we see that $\partial_{z_j}(Q)$ and $A_i\partial_{z_j}(Q)+\partial_{z_j}(H_i)$ share $0$ CM, where $i,j\in\{1,2,\ldots,m\}$ such that $i\neq j$. Consequently $\partial_{z_j}(Q)$ and $\partial_{z_j}(H_i)$ share $0$ CM, where $i,j\in\{1,2,\ldots,m\}$ such that $i\neq j$. Therefore, we can write
\bea\label{n.11} C_j\partial_{z_j}(Q)\equiv A_i\partial_{z_j}(Q)+\partial_{z_j}(H_i)\eea
and
\bea\label{n.11a} D_j\partial_{z_j}(Q)\equiv \partial_{z_j}(H_i),\eea
where $C_j, D\in\mathbb{C}\backslash \{0\}$ and $i,j\in\{1,2,\ldots,m\}$ such that $i\neq j$.\par

\medskip
First suppose $C_j=A_i$ for all $i,j\in\{1,2,\ldots,m\}$ such that $i\neq j$. Then from (\ref{n.11}), we get $\partial_{z_j}(H_i)\equiv 0$, where $i,j\in\{1,2,\ldots,m\}$ such that $i\neq j$. Clearly $H_i(z)$ is a constant, say $D_i$, where $i\in\{1,2,\ldots,m\}$. Also from (\ref{n.11a}), we get $\partial_{z_j}(Q)\equiv 0$, where $i\in\{1,2,\ldots,m\}$. This shows that $Q(z)$ is a polynomial in $\mathbb{C}^m$ of degree one. Therefore from (\ref{n.8}), we have
\bea\label{n.11b} f=e^{A_iQ+D_i},\eea
where $i=1,2,\ldots,m$. Now from (\ref{n.11b}), we get
\bea\label{n.11c} e^{(A_i-A_j)Q+D_i-D_j}=1,\eea
where $i,j\in\{1,2,\ldots,m\}$ such that $i\neq j$. Therefore from (\ref{n.11c}), one can easily deduce that $A_1=A_2=\ldots=A_m$ and $D_1=D_2=\ldots=D_m$. Consequently, we can take $f=e^{AQ+B}$, where $A\in\mathbb{C}\backslash \{0\}$ and $B\in\mathbb{C}$. Finally, we have
\[f(z)=e^{A(a_1z_1+a_2z_2+\ldots+a_mz_m)+a}\;\;\text{and}\;\;g(z)=e^{a_1z_1+a_2z_2+\ldots+a_mz_m+b},\]
where $a,A,a_1,a_2,\ldots,a_m,b\in\mathbb{C}$ such that $A\neq 0$ and $(a_1,a_2,\ldots,a_m)\neq (0,0,\ldots,0)$.\par

\medskip
Next suppose $C_j\neq A_i$ for at least one set $i,j\in\{1,2,\ldots,m\}$ such that $i\neq j$. Without loss of generality, we may assume that $C_j\neq A_i$ for all $i,j\in\{1,2,\ldots,m\}$ such that $i\neq j$. Then from (\ref{n.11}), we get
\bea\label{n.11d} (C_j-A_i)\partial_{z_j}(Q)=\partial_{z_j}(H_i),\eea
where $i,j\in\{1,2,\ldots,m\}$ such that $i\neq j$. Note that $\partial_{z_j}(H_i)$ is independent of $z_i$, where $i,j\in\{1,2,\ldots,m\}$ such that $i\neq j$. Therefore from (\ref{n.11d}), we can say that $\partial_{z_j}(Q)$ is independent of $z_j$, where $j\in\{1,2,\ldots,m\}$. Since $Q$ is a polynomial in $\mathbb{C}^m$, one can easily deduce that $Q$ is a polynomial of degree one. Therefore from (\ref{n.11d}), it follows that $H_i(z)=H(z_1,z_2,\ldots,z_{i-1},z_{i+1},\ldots,z_m)$ is also a polynomial of degree one, where $i=1,2,\ldots,m$. Clearly from (\ref{n.11}), we conclude that $\deg(A_iQ+H_i)=1$, where $i=1,2,\ldots,m$.
Finally, we can take
\[f(z)=e^{a_1z_1+a_2z_2+\ldots+a_m z_m+a}\;\;\text{and}\;\;g(z)=e^{b_1z_1+b_2z_2+\ldots+b_mz_m+b},\]
where $a,a_1,\ldots,a_m,b_1,\ldots,b_m,b\in\mathbb{C}$ such that $(a_1,\ldots,a_m)\neq (0,\ldots,0)$ and $(b_1,\ldots,b_m)\neq (0,\ldots,0)$.\par

Hence the proof is complete.
\end{proof}

\begin{lem}\label{L.3} Let $f$ and $g$ be two non-constant finite order meromorphic function in $\mathbb{C}^m$ such that $f$ and $g$ share $0$ and $\infty$ CM. If $\Phi_{0,i}\equiv 0$ for all $i=1,2,\ldots,m$, where $\Phi_{0,i}$ is defined by
\bea\label{a.1} \Phi_{0,i}=\frac{\partial_{z_i}(f)}{f}-\frac{\partial_{z_i}(g)}{g},\eea
then $f=cg$, where $c\in\mathbb{C}\backslash \{0\}$.
\end{lem}
\begin{proof} By the given conditions, we know that $f$ and $g$ share $0$ and $\infty$ CM and so $f=ge^P$, where $P$ is an entire function in $\mathbb{C}^m$. Since $\rho(f)<+\infty$ and $\rho(g)<+\infty$, using Lemma \ref{L.0}, we can say that $P$ is a polynomial in $\mathbb{C}^m$. Now using $f=ge^P$ to (\ref{a.1}), we have $\Phi_{0,i}=\partial_{z_i}(P)$. Since  $\Phi_{0,i}\equiv 0$ for all $i=1,2,\ldots,m$, we get $\partial_{z_i}(P)\equiv 0$ for all $i=1,2,\ldots,m$ and so $P$ is a constant. Therefore $f=cg$, where $c\in\mathbb{C}\backslash \{0\}$.
\end{proof}

\begin{lem}\label{L.4} Let $f$ and $g$ be two non-constant finite order meromorphic function in $\mathbb{C}^m$ such that $f$ and $g$ share $0$ and $\infty$ CM.
If $\Phi_{1,ij}\equiv 0$ for all $i,j=1,2,\ldots,m$, where $\Phi_{1,ij}$ is defined by
\bea\label{a.2} \Phi_{1,ij}=\frac{\partial^2_{z_jz_i}(f)}{\partial_{z_i}(f)}-\frac{\partial^2_{z_jz_i}(g)}{\partial_{z_i}(g)},\eea
then either $f=cg$, where $c\in\mathbb{C}\backslash \{0\}$ or
\[f=\frac{d}{1-ce^{-P}}\;\;\text{and}\;\;g=\frac{d}{e^P-c},\]
where $c,d\in\mathbb{C}\backslash \{0\}$ and $P$ is a non-constant polynomial in $\mathbb{C}^m$.
\end{lem}
\begin{proof} By the given conditions, we know that $f$ and $g$ share $0$ and $\infty$ CM and so $f=ge^P$, where $P$ is a polynomial in $\mathbb{C}^m$. If $P$ is a constant, then we have $f=cg$, where $c\in\mathbb{C}\backslash \{0\}$. Hence we assume that $P$ is a non-constant polynomial. Note that $\Phi_{1,ij}\equiv 0$ for all $i=1,2,\ldots,m$. Then from (\ref{a.2}), we get
\[\partial^2_{z_jz_i}(f)\partial_{z_i}(g)-\partial^2_{z_jz_i}(g)\partial_{z_i}(f)\equiv 0.\]

Let $F_i=\partial_{z_i}(f)$ and $G_i=\partial_{z_i}(g)$. If possible, suppose $F$ and $G$ are linearly independent. Then by Corollary 1.40 \cite{HLY1}, there exists $j\in \mathbb{Z}[1, m]$ such that
\[W(F_i,G_i)=\left|\begin{array}{ll} F_i&G_i\\\partial_{z_j}(F_i)&\partial_{z_j}(G_i)\end{array}\right|\not\equiv 0\]
and so $\partial^2_{z_jz_i}(f)\partial_{z_i}(g)-\partial^2_{z_jz_i}(g)\partial_{z_i}(f)\not\equiv 0$, which is absurd. Hence $F_i$ and $G_i$ are linearly dependent and so
\bea\label{aa.1}\partial_{z_i}(f)=c_i\partial_{z_i}(g),\eea
where $c_i\in\mathbb{C}\backslash \{0\}$ and $i=1,2,\ldots,m$. Since $f=ge^P$, we have
\bea\label{aa.2} \partial_{z_i}(f)=(\partial_{z_i}(g)+g\partial_{z_i}(P))e^P,\eea
where $i=1,2,\ldots,m$. Now using (\ref{aa.1}) and (\ref{aa.2}), we have
\bea\label{aa.3} \frac{\partial_{z_i}(g)}{g}=\frac{\partial_{z_i}(P)e^P}{c_i-e^P}\;\;\text{and}\;\;\frac{\partial_{z_i}(f)}{f}=\frac{c_i\partial_{z_i}(P)e^{-P}}{c_ie^{-P}-1},\eea
where $i=1,2,\ldots,m$. On integration, we get
\bea\label{aa.4} g(z)=\frac{e^{R_{1i}(z)}}{c_i-e^{P(z)}}\;\;\text{and}\;\;f(z)=\frac{e^{R_{2i}(z)}}{c_ie^{-P}-1},\eea
where $i=1,2,\ldots,m$ and $R_{ji}(z)=R_j(z_1,z_2,\ldots,z_{i-1},z_{i+1},\ldots,z_m)$ is a polynomial. It is easy to get from (\ref{aa.4}) that $c_1=c_2=\ldots=c_m=c$ and so
\bea\label{aa.5} f(z)=\frac{e^{R_{2i}(z)}}{ce^{-P}-1}\;\;\text{and}\;\;g(z)=\frac{e^{R_{1i}(z)}}{c-e^{P(z)}},\eea
where $i=1,2,\ldots,m$. Now from (\ref{aa.5}), we deduce that $R_{2i}-R_{2j}$ and $R_{1i}-R_{1j}$ are constants for $i,j=1,2,\ldots,m$ and $i\neq j$. Clearly $R_{11}-R_{12}$, $R_{11}-R_{13}$, $\ldots$, $R_{11}-R_{1m}$ are constants. So $\partial_{z_2}(R_{11})\equiv 0$, $\partial_{z_3}(R_{11})\equiv 0$, $\ldots$, $\partial_{z_m}(R_{11})\equiv 0$. Consequently $R_{11}$ is a constant. Similarly $R_{12}, R_{13},\ldots,R_{1m}$ are also constants. In the like manner, we can prove that $R_{21}, R_{12},\ldots,R_{2m}$ are constants.
Finally from (\ref{aa.5}), we may assume that
\beas f(z)=\frac{a}{ce^{-P}-1}\;\;\text{and}\;\;g(z)=\frac{b}{c-e^{P(z)}},\eeas
where $a,b\in\mathbb{C}\backslash \{0\}$. Hence the proof is complete.
\end{proof}

\section {{\bf Proof of Theorem \ref{t2.1}}}
\begin{proof} By the given conditions, we have
\bea\label{al.1} f=ge^P,\quad \partial_{z_i}(f)=\partial_{z_i}(g)e^{Q_i}\quad \text{and}\quad \partial^n_{z_i}(f)=\partial^n_{z_i}(g)e^{R_i},\eea
where $P$, $Q_i\;(i=1,2,\ldots,m)$ and $R_i\;(i=1,2,\ldots,m)$ are polynomials in $\mathbb{C}^m$. Here we consider the auxiliary functions $\Phi_{0,i}$ and $\Phi_{1,ij}$ defined as in (\ref{a.1}) and (\ref{a.2}) respectively.

\medskip
If $\Phi_{0,i}\equiv 0$ for $i=1,2,\ldots,m$, then by Lemma \ref{L.3}, we have $f=cg$, where $c\in\mathbb{C}\backslash \{0\}$. Henceforth we assume that $\Phi_{0,i}\not\equiv 0$ for at least one $i=1,2,\ldots,m$.

\medskip
If $\Phi_{1,ij}\equiv 0$ for $i,j=1,2,\ldots,m$, then by Lemma \ref{L.4}, we have either $f=cg$, where $c\in\mathbb{C}\backslash \{0\}$ or
\[f=\frac{d}{1-ce^{-P}}\;\;\text{and}\;\;g=\frac{d}{e^P-c},\]
where $c,d\in\mathbb{C}\backslash \{0\}$ and $P$ is a non-constant polynomial in $\mathbb{C}^m$. Henceforth, we assume that $\Phi_{1,ij}\not\equiv 0$ for at least one $(i,j)\in\{i=1,2,\ldots,m\}\times \{1,2,\ldots,m\}$.

\medskip
Define $I=I_f\cup I_g$ and $S=\sideset{}{_{a\in\{0,\infty\}}}{\bigcup}(\operatorname{supp}\mu_f^a)_s\cup(\operatorname{supp}\mu_g^a)_s$,
where $A_s$ denotes the set of singular points of the analytic set $A$. Note that $\dim (I\cup S)\leq m-2$. Take $z_0\in\mathbb{C}^m-(I\cup S)$.

\smallskip
Suppose $p=\mu^0_f(z_0)>0$. By the given condition, we have $p=\mu^0_g(z_0)>0$.
Note that $f^{-1}(\{0\})_s$ is an analytic subset of $\mathbb{C}^m$ of dimension $\leq m-2$, where $f^{-1}(\{0\})_s$ is the set of singular points of $f^{-1}(\{0\})$. Clearly $z_0\in \mathbb{C}^m-f^{-1}(\{0\})_s$. Then there is a holomorphic coordinate system $\left(U ; \varphi_1, \ldots, \varphi_m\right)$ of $z_0$ in $\mathbb{C}^m$ such that
\[U \cap f^{-1}(\{0\})=\left\{z \in U \mid \varphi_1(z)=0\right\}\;\;\text{and}\;\;\varphi_j(z_0)=0\]
for $j=1,2,\ldots,m$ (see Lemma 2.3 \cite{BCL} and the proof of Lemma 2.3 \cite{FL1}). Therefore biholomorphic coordinate transformation $z_j=z_j(\varphi_1, \ldots, \varphi_m)$, $j=1, \ldots, m$
near $0$ exists such that $z_0=z(0)=(z_1(0), \ldots, z_m(0))$.
So we can write
\bea\label{al.5} f=a_p\varphi_1^{p}+a_{p+1}\varphi_1^{p+1}+a_{p+2}\varphi_1^{p+2}+\varphi_1^{p+3}\hat{f}\left(\varphi_1, \ldots, \varphi_m\right)\eea
and
\bea\label{al.6} g=b_p\varphi_1^{p}+b_{p+1}\varphi_1^{p+1}+b_{p+2}\varphi_1^{p+2}+\varphi_1^{p+3}\hat{g}\left(\varphi_1, \ldots, \varphi_m\right),\eea
where $a_p(\neq 0), a_{p+1}, a_{p+2}, b_p(\neq 0), b_{p+1},b_{p+2}\in\mathbb{C}$ and $\hat{f}$ and $\hat g$ are holomorphic function near $0$.
Clearly (\ref{al.5}) and (\ref{al.6}) give
\bea\label{al.5a}\partial_{z_i}(f)=pa_p\varphi_1^{p-1}\partial_{z_i}(\varphi_1)+(p+1)a_{p+1}\varphi_1^{p}\partial_{z_i}(\varphi_1)+(p+2)a_{p+2}\varphi_1^{p+1}\partial_{z_i}(\varphi_1)+\ldots,\eea
\bea\label{al.5b}\partial^2_{z_jz_i}(f)&=&p(p-1)a_p\varphi_1^{p-2}\partial_{z_j}(\varphi_1)\partial_{z_i}(\varphi_1)+pa_p\varphi_1^{p-1}\partial^2_{z_jz_i}(\varphi_1)\\&&+(p+1)pa_{p+1}\varphi_1^{p-1}\partial_{z_j}(\varphi_1)\partial_{z_i}(\varphi_1)+(p+1)a_{p+1}\varphi_1^{p}\partial^2_{z_jz_i}(\varphi_1)\nonumber\\&&+(p+2)(p+1)a_{p+2}\varphi_1^{p}\partial_{z_j}(\varphi_1)\partial_{z_i}(\varphi_1)+(p+2)a_{p+2}\varphi_1^{p+1}\partial^2_{z_jz_i}(\varphi_1)+\ldots,\nonumber\eea
\beas\partial^2_{z_i}(f)&=&p(p-1)a_p\varphi_1^{p-2}(\partial_{z_i}(\varphi_1))^2+pa_p\varphi_1^{p-1}\partial^2_{z_i}(\varphi_1)+
(p+1)pa_{p+1}\varphi_1^{p-1}(\partial_{z_i}(\varphi_1))^2\\&&+(p+1)a_{p+1}\varphi_1^{p}\partial^2_{z_i}(\varphi_1)+(p+2)(p+1)a_{p+2}\varphi_1^{p}(\partial_{z_i}(\varphi_1))^2\nonumber\\&&+(p+2)a_{p+2}\varphi_1^{p+1}\partial^2_{z_i}(\varphi_1)+\ldots,\nonumber\eeas
\bea\label{al.5c}\partial^n_{z_i}(f)&=&p(p-1)\ldots(p-n+1)a_p\varphi_1^{p-n}(\partial_{z_i}(\varphi_1))^n+\ldots\\&&
(p+1)p\ldots(p-n+2)a_{p+1}\varphi_1^{p-n+1}(\partial_{z_i}(\varphi_1))^n\nonumber\\&&+(p+2)(p+1)\ldots(p-n+3)a_{p+2}\varphi_1^{p-n+2}(\partial_{z_i}(\varphi_1))^n+\ldots,\nonumber\eea
\bea\label{al.6a}\partial_{z_i}(g)=pb_p\varphi_1^{p-1}\partial_{z_i}(\varphi_1)+(p+1)b_{p+1}\varphi_1^{p}\partial_{z_i}(\varphi_1)+(p+2)b_{p+2}\varphi_1^{p+1}\partial_{z_i}(\varphi_1)+\ldots,\eea
\bea\label{al.6b}\partial^2_{z_jz_i}(g)&=&p(p-1)b_p\varphi_1^{p-2}\partial_{z_j}(\varphi_1)\partial_{z_i}(\varphi_1)+pb_p\varphi_1^{p-1}\partial^2_{z_jz_i}(\varphi_1)\\&&+(p+1)pb_{p+1}\varphi_1^{p-1}\partial_{z_i}(\varphi_1)\partial_{z_i}(\varphi_1)+(p+1)b_{p+1}\varphi_1^{p}\partial^2_{z_jz_i}(\varphi_1)\nonumber\\&&+(p+2)(p+1)b_{p+2}\varphi_1^{p}\partial_{z_j}(\varphi_1)\partial_{z_i}(\varphi_1)+(p+2)b_{p+2}\varphi_1^{p+1}\partial^2_{z_jz_i}(\varphi_1)+\ldots,\nonumber\eea
and
\bea\label{al.6c}\partial^n_{z_i}(g)&=&p(p-1)\ldots(p-n+1)b_p\varphi_1^{p-n}(\partial_{z_i}(\varphi_1))^n+\ldots\\&&
(p+1)p\ldots(p-n+2)b_{p+1}\varphi_1^{p-n+1}(\partial_{z_i}(\varphi_1))^n\nonumber\\&&+(p+2)(p+1)\ldots(p-n+3)b_{p+2}\varphi_1^{p-n+2}(\partial_{z_i}(\varphi_1))^n+\ldots.\nonumber\eea

Now we consider following two cases.\par

\medskip
{\bf Case 1.} Suppose $\partial_{z_i}(f)/f$ is a rational function for all $i=1,2,\ldots,m$. Then by Lemma \ref{L.2}, we get the desired conclusions.\par

\medskip
{\bf Case 2.} Suppose $\partial_{z_i}(f)/f$ is transcendental for atleast one $i\in\{1,2,\ldots,m\}$. Without loss of generality, we may assume that $\partial_{z_i}(f)/f$ is transcendental for all $i=1,2,\ldots,m$. By Lemma \ref{l2.1}, we have
\[m\left(r,\frac{\partial_{z_i}(f)}{f}\right)=O(\log r).\]

Also it is easy to deduce that
\[N\left(r,\frac{\partial_{z_i}(f)}{f}\right)\leq \ol N(r,f)+\ol N(r,0;f)=\ol N(r,f)+\sideset{}{_{k=1}^{n-1}}{\sum}\ol N_{=k}(r,0;f)+\ol N_{(n}(r,0;f).\]

Consequently
\bea\label{al.2}T\left(r,\frac{\partial_{z_i}(f)}{f}\right)\leq \ol N(r,f)+\sideset{}{_{k=1}^{n-1}}{\sum}\ol N_{=i}(r,0;f)+\ol N_{(n}(r,0;f)+O(\log r).\eea

\smallskip
Let
\[S_k=\{z\in\mathbb{C}^m:\ol{\mu}^0_{f=k}(z)>0\}, S_{(n}=\{z\in\mathbb{C}^m:\ol{\mu}^0_{f(n}(z)>0\}\;\text{and}\; S=\{z\in\mathbb{C}^m:\ol{\mu}^{\infty}_{f}(z)>0\},\]
where $k=1,2,\ldots,n-1$. Since $\partial_{z_i}(f)/f$ is transcendental, from (\ref{al.2}) it follows that the set $\bigcup_{k=1}^{n-1} S_k\cup S_{(n}\cup S$ is not algebraic.

We now consider following two sub-cases.\par

\medskip
{\bf Sub-case 2.1.} Suppose the set $S_k$ is not algebraic for at least one value of $k\in\{1,2,\ldots,n-1\}$.
Let
\bea\label{al.3} \Psi_{k,ij}=(k+1)\Phi_{0,j}-k\Phi_{1,ij}=(k+1)\left(\frac{\partial_{z_j}(f)}{f}-\frac{\partial_{z_j}(g)}{g}\right)-k\left(\frac{\partial^2_{z_jz_i}(f)}{\partial_{z_i}(f)}-\frac{\partial^2_{z_jz_i}(g)}{\partial_{z_i}(g)}\right),\eea
$k=1,2,\ldots,n-1$ and $i,j=1,2,\ldots,m$. On the other hand from (\ref{al.1}) and (\ref{al.3}), we deduce that
\bea\label{al.4} \Psi_{k,ij}=(k+1)\partial_{z_j}(P)-k\partial_{z_j}(Q_i)\eea
$k=1,2,\ldots,n-1$ and $i,j=1,2,\ldots,m$.

Now we divide following two sub-cases.\par

\medskip
{\bf Sub-case 2.1.1.} Suppose $\Psi_{k,ij}\not\equiv 0$ for at least one $(k,i,j)\in \{1,2,\ldots,n-1\}\times \{1,2,\ldots,m\}\times\{1,2,\ldots,m\}$. Let $z_0\in S_k$ such that $k=\mu^0_f(z_0)>0$. By the given condition, we have $k=\mu^0_g(z_0)>0$. Let us take $p=k$. Now we apply (\ref{al.5})-(\ref{al.5b}), (\ref{al.6a}) and (\ref{al.6b}) to (\ref{al.3}) and then by a elementary calculation, we can conclude that $\Psi_{k,ij}(z_0)=0$. Consequently the set $S_k$ is algebraic, which is absurd.\par

\medskip
{\bf Sub-case 2.1.2.} Suppose $\Psi_{k,ij}\equiv 0$ for all $(k,i,j)\in \{1,2,\ldots,n-1\}\times \{1,2,\ldots,m\}\times\{1,2,\ldots,m\}$. Then from (\ref{al.4}), we get
\bea\label{al.8a} (k+1)\partial_{z_j}(P)-k\partial_{z_j}(Q_i)\equiv 0,\eea
where $k=1,2,\ldots,n-1$ and $i,j=1,2,\ldots,m$. On integration, we get
\bea\label{al.9} Q_i(z)=\frac{k+1}{k} P(z)+\tilde R_j(z),\eea
where $\tilde R_j(z)=\tilde R(z_1,z_2,\ldots,z_{j-1}, z_{j+1},\ldots,z_m)$ is a polynomial. Clearly from (\ref{al.9}), we get
\[Q_1(z)=\frac{k+1}{k} P(z)+\tilde R_1(z), Q_1(z)=\frac{k+1}{k} P(z)+\tilde R_2(z),\ldots,Q_1(z)=\frac{k+1}{k} P(z)+\tilde R_m(z)\]
which demands that $\tilde R_1(z)=\tilde R_2(z)=\ldots=\tilde R_m(z)=A\in\mathbb{C}$. Therefore from (\ref{al.9}), we get
\bea\label{al.9a} Q_i(z)=\frac{k+1}{k} P(z)+A.\eea

Now from (\ref{al.1}), we have
\bea\label{al.10} \frac{\partial_{z_i}(f)}{f}=\frac{\partial_{z_i}(g)}{g}+\partial_{z_i}(P)\eea
where $i=1,2,\ldots,m$. Again from (\ref{al.1}), we get
\bea\label{al.11} \frac{\partial_{z_i}(f)}{f}=\frac{\partial_{z_i}(g)}{g}e^{Q_i-P},\eea
where $i=1,2,\ldots,m$. Therefore from (\ref{al.10}) and (\ref{al.11}), we obtain
\bea\label{al.12} \frac{\partial_{z_i}(g)}{g}=\frac{-\partial_{z_i}(P)}{1-e^{Q_i-P}}\;\;\text{and}\;\;\frac{\partial_{z_i}(f)}{f}=\frac{-\partial_{z_i}(P)}{e^{-Q_i+P}-1},\eea
where $i=1,2,\ldots,m$.

Note that the set $S_k$ is not algebraic for at least one value of $k\in\{1,2,\ldots,n-1\}$. Since $\partial_{z_i}(g)/g$ has only simple poles, it follows from (\ref{al.12}) that $\deg(Q_i-P)=1$. Consequently from (\ref{al.8a}), we can conclude that $\deg(P)=\deg(Q_i)=1$. Again using (\ref{al.9a}) to (\ref{al.12}), we get
\bea\label{al.13} \frac{\partial_{z_i}(g)}{g}=k\frac{\frac{1}{k}\partial_{z_i}(P)e^{-\frac{1}{k}P-A}}{1-e^{-\frac{1}{k}P-A}}\;\;\text{and}\;\;\frac{\partial_{z_i}(f)}{f}=k\frac{\frac{1}{k}\partial_{z_i}(P)e^{\frac{1}{k}P+A}}{e^{\frac{1}{k}P+A}-1},\eea
where $i,j=1,2,\ldots,m$. On integration, we get from (\ref{al.13}) that
\bea\label{al.14} g(z)=\left(1-e^{-\frac{1}{k}P(z)-A}\right)^ke^{\hat R_{1i}(z)}\;\;\text{and}\;\;f(z)=\left(e^{\frac{1}{k}P(z)+A}-1\right)^ke^{\hat R_{2i}(z)},\eea
where $i=1,2,\ldots,m$ and $\hat R_{ji}(z)=\hat R_j(z_1,z_2,\ldots,z_{i-1},z_{i+1},\ldots,z_m)$ is a polynomial. Now proceeding in the same way as done in the proof of Lemma \ref{L.4}, we can prove that $\hat R_{1i}(z)$ and $\hat R_{2i}(z)$ are constants.
Consequently from (\ref{al.14}), we have
\bea\label{al.16} g(z)=A_i\left(1-e^{-\frac{1}{k}P(z)-A}\right)^k\;\;\text{and}\;\;f(z)=B_i\left(e^{\frac{1}{k}P(z)+A}-1\right)^k,\eea
where $A_i,B_i\in\mathbb{C}\backslash \{0\}$ and $i=1,2,\ldots,m$. Again from (\ref{al.16}), one can easily prove that $A_1=A_2=\ldots=A_m$ and $B_1=B_2=\ldots=B_m$. Let us take
\[P(z)=a_1z_1+a_2z_2+\ldots+a_mz_m+a_0,\]
where $a_0,a_1,a_2,\ldots,a_m\in\mathbb{C}$ such that $(a_1,a_2,\ldots,a_m)\neq (0,0,\ldots,0)$. Now in view of (\ref{al.16}), we may assume that
\bea\label{al.17a} f(z)=a\left(1-be^{\frac{1}{k}(a_1z_1+a_2z_2+\ldots+a_mz_m)}\right)^k\;\;\text{and}\;\;g(z)=c\left(e^{-\frac{1}{k}(a_1z_1+a_2z_2+\ldots+a_mz_m)}-b\right)^k,\eea
where $a,b,c\in\mathbb{C}\backslash \{0\}$. Applying the binomial theorem, we get
\bea\label{al.17} f(z)=a\left\lbrace 1-kbe^{\frac{1}{k}(a_1z_1+\ldots+a_mz_m)}+\binom{k}{2}b^2e^{\frac{2}{k}(a_1z_1+\ldots+a_mz_m)}-\ldots+(-b)^k e^{a_1z_1+\ldots+a_mz_m}\right\rbrace\eea
and
\bea\label{al.18} g(z)&=&c\left\lbrace e^{-(a_1z_1+\ldots+a_mz_m)}-kbe^{-\frac{k-1}{k}(a_1z_1+\ldots+a_mz_m)}+\binom{k}{2}b^2e^{-\frac{k-2}{k}(a_1z_1+\ldots+a_mz_m)}\right.\\&&\left.-\ldots+k(-b)^{k-1}e^{-\frac{1}{k}(a_1z_1+\ldots+a_mz_m)}+(-b)^k\right\rbrace\nonumber.\eea

Now differentiating partially (\ref{al.17}) and (\ref{al.18}) $n$-times with respect to $z_i$, we get respectively
\bea\label{al.19} \partial^n_{z_i}(f(z))&=&a(a_i/k)^n\Big(-kb e^{\frac{1}{k}(a_1z_1+\ldots+a_mz_m)}+\binom{k}{2}b^2 2^ne^{\frac{2}{k}(a_1z_1+\ldots+a_mz_m)}\\&&-\ldots+(-b)^k k^n e^{a_1z_1+\ldots+a_mz_m}\Big)\nonumber\eea
and
\bea\label{al.20} \partial^n_{z_i}(g(z))&=&c(-a_i/k)^n\Big(k^ne^{-(a_1z_1+\ldots+a_mz_m)}-kb(k-1)^ne^{-\frac{k-1}{k}(a_1z_1+\ldots+a_mz_m)}\\&&+\binom{k}{2}b^2(k-2)^ne^{-\frac{k-2}{k}(a_1z_1+\ldots+a_mz_m)}-\ldots+k(-b)^{k-1}e^{-\frac{1}{k}(a_1z_1+\ldots+a_mz_m)}\Big)\nonumber.\eea

Since $\partial^n_{z_i}(f(z))$ and $\partial^n_{z_i}(g(z))$ share $0$ CM, from (\ref{al.19}) and (\ref{al.20}), one can easily conclude that $k=1$. Therefore from (\ref{al.17a}), we have
\beas f(z)=a\left(1-be^{(a_1z_1+a_2z_2+\ldots+a_mz_m)}\right)\;\;\text{and}\;\;g(z)=c\left(e^{-(a_1z_1+a_2z_2+\ldots+a_mz_m)}-b\right),\eeas
where $a,b,c\in\mathbb{C}\backslash \{0\}$ and $a_1,a_2,\ldots,a_m\in\mathbb{C}$ such that $(a_1,a_2,\ldots,a_m)\neq (0,0,\ldots,0)$.\par

\medskip
{\bf Sub-case 2.2.} Suppose the set $S_k$ is algebraic for all $k\in\{1,2,\ldots,n-1\}$. In this case, the set $S_{(n}\cup S$ is not algebraic. We define
\bea\label{a.3} \Phi_{n,ij}=\frac{\partial^{n+1}_{z_jz_i^n}(f)}{\partial^n_{z_i}(f)}-\frac{\partial^{n+1}_{z_jz_i^n}(g)}{\partial^n_{z_i}(g)},\eea
where $i=1,2,\ldots,m$ and $j=1,2,\ldots,m$.

\smallskip
Suppose $p=\mu^{0}_f(z_0)\geq n$. By the given condition, we have $p=\mu^{0}_g(z_0)\geq n$. Set
\bea\label{al.24a} A_0:=A_0(z_0)=\frac{a_{p+1}}{a_p}-\frac{b_{p+1}}{b_p}\eea
\bea\label{al.24b} B_0:=B_0(z_0)=\frac{a_{p+2}}{a_p}-\frac{b_{p+2}}{b_p}\eea
and
\bea\label{al.24c} C_0:=C_0(z_0)=A_0(z_0)\left(\frac{a_{p+1}}{a_p}+\frac{b_{p+1}}{b_p}\right).\eea

Now using (\ref{al.5})-(\ref{al.6c}) to (\ref{a.1}), (\ref{a.2}) and (\ref{a.3}), we get respectively
\bea\label{al.25a} \Phi_{0,j}(z)=A_0+2B_0\varphi_1-C_0 \varphi_1+O(\varphi_1^2),\eea
\bea\label{al.25b} \Phi_{1,ij}(z)=A_0\frac{p+1}{p}+2B_0\frac{p+2}{p}\varphi_1-C_0\frac{(p+1)^2}{p^2} \varphi_1+O(\varphi_1^2)\eea
and
\bea\label{al.25c} \Phi_{n,ij}(z)=A_0 \frac{\binom{p+1}{n}}{\binom{p}{n}}+2B_0\frac{\binom{p+2}{n}}{\binom{p}{n}}\varphi_1-C_0 \frac{\binom{p+1}{n}^2}{\binom{p}{n}^2} \varphi_1+O(\varphi_1^2). \eea

\smallskip
Suppose $p=\mu^{\infty}_f(z_0)>0$. By the given condition, we have $p=\mu^{\infty}_g(z_0)>0$.
Then there is a holomorphic coordinate system $\left(U; u_1, \ldots, u_m\right)$ of $z_0$ in $\mathbb{C}^m-(I\cup S)$ such that
\[U \cap \text{supp} \mu^{\infty}_f=\left\{z \in U \mid u_1(z)=0\right\}\;\;\text{and}\;\;u_j(z_0)=0\]
for $j=1,2,\ldots,m$. Therefore biholomorphic coordinate transformation $z_j=z_j(u_1, \ldots, u_m)$, $j=1, \ldots, m$
near $0$ exists such that $z_0=z(0)=(z_1(0), \ldots, z_m(0))$. Therefore, we can write
\bea\label{al.27} f=c_p\varphi_1^{-p}+c_{p+1}\varphi_1^{-p+1}+c_{p+2}\varphi_1^{-p+2}+\varphi_1^{-p+3}\tilde{f}\left(u_1, \ldots,u_m\right)\eea
and
\bea\label{al.28} g=d_p\varphi_1^{-p}+d_{-p+1}\varphi_1^{p+1}+d_{-p+2}\varphi_1^{-p+2}+\varphi_1^{-p+3}\tilde {g}\left(u_1, \ldots,u_m\right),\eea
where $c_p(\neq 0), c_{p+1}, c_{p+2}, d_p(\neq 0), d_{p+1},d_{p+2}\in\mathbb{C}$ and $\tilde {f}$ and $\tilde g$ are holomorphic function near $0$.

Set
\bea\label{al.28a} A_{\infty}:=A_{\infty}(z_0)=\frac{c_{p+1}}{c_p}-\frac{d_{p+1}}{d_p}\eea
\bea\label{al.28b} B_{\infty}:=B_{\infty}(z_0)=\frac{c_{p+2}}{c_p}-\frac{d_{p+2}}{d_p}\eea
and
\bea\label{al.28c} C_{\infty}:=C_{\infty}(z_0)=A_{\infty}(z_0)\left(\frac{c_{p+1}}{c_p}+\frac{d_{p+1}}{d_p}\right).\eea

In this case also, we have
\bea\label{al.29a} \Phi_{0,j}(z)=A_{\infty}+2B_{\infty}u_1-C_{\infty} u_1+O(u_1^2),\eea
\bea\label{al.29b} \Phi_{1,ij}(z)=A_{\infty}\frac{p+1}{p}+2B_{\infty}\frac{p+2}{p}u_1-C_{\infty}\frac{(p+1)^2}{p^2} u_1+O(u_1^2)\eea
and
\bea\label{al.29c} \Phi_{n,ij}(z)=A_{\infty} \frac{\binom{p+1}{n}}{\binom{p}{n}}+2B_{\infty}\frac{\binom{p+2}{n}}{\binom{p}{n}}u_1-C_{\infty} \frac{\binom{p+1}{n}^2}{\binom{p}{n}^2}u_1+O(u_1^2).\eea

\smallskip
Let
\[\hat S_0=\{z\in\mathbb{C}^m: \ol{\mu}^0_{f(n}(z)>0\;\;\text{and}\;\;A_0(z)=0\}\;\;\text{and}\;\;\hat S_{\infty}=\{z\in\mathbb{C}^m: \ol{\mu}^{\infty}_{f,1}(z)\geq 1\;\;\text{and}\;\;A_{\infty}(z)=0\}.\]

Now we consider following two sub-cases.\par

\medskip
{\bf Sub-case 2.2.1.} Suppose the set $\hat S_0\cup \hat S_{\infty}$ is not algebraic. By assumption, we have $\Phi_{0,i}\not\equiv 0$ for at least one value of $i\in\{1,2,\ldots,m\}$. Again using (\ref{al.1}) to (\ref{a.1}), we get $\Phi_{0,i}=\partial_{z_i}(P)$ and so  $T(r,\Phi_{0,i})=O(\log r)$ for $i=1,2,\ldots,m$. Now using (\ref{al.24a}) to (\ref{al.25a}) and (\ref{al.28a}) to (\ref{al.29a}), we get
\[\ol N_{(n}(r,0;f)+\ol N(r,f)\leq N(r,0;\Phi_{0,i})\leq T(r,\Phi_{0,i})+O(1)=O(\log r),\]
which shows that the set $S_{(n}\cup S$ is algebraic. Therefore, we get a contradiction.

\smallskip
{\bf Sub-case 2.2.2.} Suppose the set $\hat S_0\cup \hat S_{\infty}$ is algebraic. In this case, the set $\tilde S_0\cup \tilde S_{\infty}$ is not algebraic, where
\[\tilde S_0=\{z\in\mathbb{C}^m: \ol{\mu}^0_{f(n}(z)>0\;\;\text{and}\;\;A_0(z)\neq 0\}\;\;\text{and}\;\; \tilde S_{\infty}=\{z\in\mathbb{C}^m: \ol{\mu}^{\infty}_{f,1}(z)> 0\;\;\text{and}\;\;A_{\infty}(z)\neq 0\}.\]

\smallskip
Also in this case, $\Phi_{0,i}\not\equiv 0$ and $\Phi_{s,ij}\not\equiv 0$ for all $i,j=1,2,\ldots,m$ and $s=1,n$. Now using (\ref{al.1}) to (\ref{a.1}), (\ref{a.2}) and (\ref{a.3}), we get respectively $\Phi_{0,i}=\partial_{z_i}(P)$, $\Phi_{1,ij}=\partial_{z_j}(Q_i)$ and $\Phi_{n,ij})=\partial_{z_j}(R_i)$. Therefore $T(r,\Phi_{0,i})=O(\log r)$ and $T(r,\Phi_{s,ij})=O(\log r)$, where $s=1,n$.

\smallskip
Let
\bea\label{la.1} H_{ij}=\frac{n-1}{\Phi_{0,j}}-\frac{n}{\Phi_{1,ij}}+\frac{1}{\Phi_{n,ij}},\eea
where $i,j=1,2,\ldots,m$. Since $T(r,\Phi_{0,i})=O(\log r)$ and $T(r,\Phi_{s,ij})=O(\log r)$, where $s=1,n$, it follows that $T(r,H_{ij})=O(\log r)$.

\smallskip
Now if $z_0\in\tilde S_0$, then using (\ref{al.25a})-(\ref{al.25c}) to (\ref{la.1}), we get $H_{ij}(z_0)=0$. Again if $z_0\in\tilde S_{\infty}$, then using (\ref{al.29a})-(\ref{al.29c}) to (\ref{la.1}), we get $H_{ij}(z_0)=0$. Since $T(r,H_{ij})=O(\log r)$, we can conclude that the set $\tilde S_0\cup \tilde S_{\infty}$ is algebraic, which is absurd. Hence $H_{ij}\equiv 0$ and so from (\ref{al.1}), we get
\bea\label{la.2} \frac{n-1}{\Phi_{0,j}}-\frac{n}{\Phi_{1,ij}}+\frac{1}{\Phi_{n,ij}}\equiv 0.\eea

\smallskip
Now if $z_0\in\tilde S_0$, then using (\ref{al.25a})-(\ref{al.25c}) to (\ref{la.2}), we get $B_0(z_0)=0$. Again if $z_0\in\tilde S_{\infty}$, then using (\ref{al.29a})-(\ref{al.29c}) to (\ref{la.2}), we get $B_{\infty}(z_0)=0$.
Therefore if $z_0\in\tilde S_0$, then from (\ref{al.25a})-(\ref{al.25c}), we get
\bea\label{la.3} \frac{1}{\Phi_{0,j}(z)}=\frac{1}{A_0}+\frac{C_0}{A_0^2}\varphi_1+O(\varphi_1^2)\eea
and
\bea\label{la.3a} \frac{1}{\phi_{s,ij}(z)}=\frac{\binom{p}{s}}{\binom{p+1}{s}A_0}+\frac{C_0}{A_0^2}\varphi_1+O(\varphi_1^2),\eea
where $s=1,n$ and if $z_0\in\tilde S_{\infty}$, then from (\ref{al.29a})-(\ref{al.29c}), we get
\bea\label{la.4} \frac{1}{\Phi_{0,j}(z)}=\frac{1}{A_{\infty}}+\frac{C_{\infty}}{A_{\infty}^2}u_1+O(u_1^2).\eea
and
\bea\label{la.4a} \frac{1}{\Phi_{s,ij}(z)}=\frac{\binom{p}{s}}{\binom{p+1}{s}A_{\infty}}+\frac{C_{\infty}}{A_{\infty}^2}u_1+O(u_1^2),\eea
where $s=1,n$. Let
\bea\label{la.5} G_{s,ij}:=\partial_{z_i}\left(\frac{1}{\Phi_{0,j}}-\frac{1}{\Phi_{s,ij}}\right),\eea
where $s=1,n$.

\smallskip
Since $T(r,\Phi_{0,i})=O(\log r)$ and $T(r,\Phi_{s,ij})=O(\log r)$, where $s=1,n$, it follows that $T(r,G_{s,ij})=O(\log r)$.
Now using (\ref{la.3})-(\ref{la.4a}) to (\ref{la.5}), we deduce that $\partial_{z_i}(G_{s,ij}(z_0))=0$ for $s=1,n$. Consequently, the set $\tilde S_0\cup \tilde S_{\infty}$ is algebraic, which is absurd. Hence $G_{s,ij}\equiv 0$ for $s=1,n$ and so
\bea\label{la.6}\frac{1}{\Phi_{0,j}}-\frac{1}{\Phi_{s,ij}}=\frac{s}{(p+1)A}=sM\eea
is a constant, where $M\in\mathbb{C}\backslash \{0\}$.

We know that $\infty$ is a Picard exceptional value of $\Phi_{0,j}$, $\Phi_{1,ij}$ and $\Phi_{n,ij}$. Therefore from (\ref{la.5}), it is easy to verify that $1/M$ and $1/nM$ are Picard exceptional values of $\Phi_{0,j}$. Consequently $\Phi_{0,j}$ is a constant. Now from (\ref{la.6}), one can conclude that $\Phi_{1,ij}$ and $\Phi_{n,ij}$ are both constants.

\smallskip
On the other hand from (\ref{la.3}) and (\ref{la.4}), we can deduce that $C_0=0$ and $C_{\infty}=0$ and so by (\ref{al.24c}) and (\ref{al.28c}), we have $a_{p+1}/a_p=-b_{p+1}/b_p$ and $c_{p+1}/c_p=-d_{p+1}/d_p$. Thus
\bea\label{la.7} A_0=2a_{p+1}/a_p\;\;\text{and}\;\;A_{\infty}=2c_{p+1}/c_p.\eea

\smallskip
Let $z_0\in\tilde S_0$. The coefficients of $\varphi_1^2$ in the expansions given by (\ref{al.25a})-(\ref{al.25c}) of $\Phi_{0,i}$, $\Phi_{1,ij}$ and $\Phi_{n,ij}$ about $z_0$ respectively, yield a homogeneous linear equation system. It is easy to see that $(a_{p+1}/a_p)^3=0$. Then from (\ref{la.7}), we have $A_0(z_0)=0$, contradicting the choice of $z_0$ (see proof of Theorem 4 \cite{LK}).

\smallskip
Let $z_0\in\tilde S_{\infty}$. The coefficients of $u_1^2$ in the expansions given by (\ref{al.29a})-(\ref{al.29c}) of $\Phi_{0,i}$, $\Phi_{1,ij}$ and $\Phi_{n,ij}$ about $z_0$ respectively, yield a homogeneous linear equation system. It is easy to see that $(a_{p+1}/a_p)^3=0$. Then from (\ref{la.7}), we have $A_{\infty}(z_0)=0$, contradicting the choice of $z_0$.

This completes the proof.
\end{proof}

\vspace{0.1in}
{\bf Compliance of Ethical Standards:}\par

{\bf Conflict of Interest.} The authors declare that there is no conflict of interest regarding the publication of this paper.\par

{\bf Data availability statement.} Data sharing not applicable to this article as no data sets were generated or analysed during the current study.

{\bf Funding.} This work was supported by
the National Natural Science Foundation of China (12161074) and the Foundation
of Education Department of Jiangxi (GJJ190876, GJJ202303, GJJ201813, GJJ201343) of China, and the Talent Introduction
Research Foundation of Suqian University (106-CK00042/028), Suqian Sci \& Tech Program Grant No. M202206 and sponsored by Qing Lan Project of Jiangsu province and Suqian Talent Xiongying Plan of Suqian.

\end{document}